\documentclass[10pt]{article}
\usepackage{latexsym, amsfonts}
\usepackage{amsmath}
\usepackage{amssymb}
\usepackage{amsthm}
\usepackage{mathrsfs}
\usepackage{mathabx}
\usepackage{dsfont}
\usepackage{bbm}
\usepackage{cite}
\usepackage{titlesec}
\usepackage{blindtext}
\usepackage{xcolor}
\usepackage{tikz-cd}
\allowdisplaybreaks

\newcommand{\be}{\begin{equation}}
\newcommand{\ee}{\end{equation}}
\newcommand{\bea}{\begin{eqnarray}}
\newcommand{\eea}{\end{eqnarray}}
\newcommand{\bean}{\begin{eqnarray*}}
\newcommand{\eean}{\end{eqnarray*}}
\newcommand{\brray}{\begin{array}}
\newcommand{\erray}{\end{array}}
\newcommand{\ben}{\begin{equation}{nonumber}}
\newcommand{\een}{\end{equation}{nonumber}}

\newtheorem{dfn}{Definition}[section]
\newtheorem{thm}[dfn]{Theorem}
\newtheorem{lmma}[dfn]{Lemma}
\newtheorem{ppsn}[dfn]{Proposition}
\newtheorem{crlre}[dfn]{Corollary}
\newtheorem{xmpl}[dfn]{Example}
\newtheorem{rmrk}[dfn]{Remark}
\newcommand{\bdfn}{\begin{dfn}}
\newcommand{\bthm}{\begin{thm}}
\newcommand{\blmma}{\begin{lmma}}
\newcommand{\bppsn}{\begin{ppsn}}
\newcommand{\bcrlre}{\begin{crlre}}
\newcommand{\bxmpl}{\begin{xmpl}}
\newcommand{\brmrk}{\begin{rmrk}}
\newcommand{\bpmat}{\begin{pmatrix}}
\newcommand{\edfn}{\end{dfn}}
\newcommand{\ethm}{\end{thm}}
\newcommand{\elmma}{\end{lmma}}
\newcommand{\eppsn}{\end{ppsn}}
\newcommand{\ecrlre}{\end{crlre}}
\newcommand{\exmpl}{\end{xmpl}}
\newcommand{\ermrk}{\end{rmrk}}
\newcommand{\epmat}{\end{pmatrix}}

\newcommand{\IC}{\mathbb{C}}

\newcommand{\al}{\alpha}
\newcommand{\bta}{\beta}

\newcommand{\cla}{{\cal A}}
\newcommand{\clb}{{\cal B}}
\newcommand{\clc}{{\cal C}}
\newcommand{\cld}{{\cal D}}

\newcommand{\clf}{{\cal F}}
\newcommand{\clg}{{\cal G}}
\newcommand{\clh}{{\cal H}}

\newcommand{\clk}{{\cal K}}

\newcommand{\clq}{{\cal Q}}

\newcommand{\cls}{{\cal S}}

\newcommand{\kac}{\mathrm{kac}}

\newcommand{\M}{\mathrm{M}}

\newcommand{\C}{\mathds{C}}

\def\a*{{\cal A}_{h,*}}
\def\B{{\cal B}(h)}
\def\B1{{\cal B}_1(h)}
\def\b{{\cal B}^{\rm s.a.}(h)}
\def\b1{{\cal B}^{\rm s.a.}_1(h)}

\newcommand{\ot}{\otimes}

\newcommand{\raro}{\rightarrow}

\def \qed {$\Box$}

\begin{document}
\begin{center}
{\large {\bf  Outer action of representation category of discrete quantum groups: a case study }}\\
by\\
{\large Debashish Goswami{\footnote{ Partially supported by JC Bose National Fellowship given by SERB, Govt. of India.}}
},
{\large Suchetana Samadder{\footnote{ Support by a research fellowship from Indian Statistical Institute is gratefully acknowledged.}}
}\\
{ Stat-Math Unit, Kolkata,}\\
{ Indian Statistical Institute,}\\
{ 203, B. T. Road, Kolkata 700 108, India}\\
{e-mails: debashish\_goswami@yahoo.co.in, suchetana.smdr@gmail.com }\\ 

\end{center}

\begin{abstract}
\noindent Given any tracial von Neumann algebra $\cla$, the first author has defined a unitary tensor functor $\Phi:{\rm Rep}(\clq_{\rm aut}(\cla))\to {\rm Bimod}(\cla)$ in his paper \cite{krp_insa}. In this paper, we consider $\cla$ to be the underlying von Neumann algebra of $C^*(S_3)$ and show that given any DQG $\clq$, which contains $C^*(S_3)$ as a quantum subgroup and has an outer action $\tilde{\Phi}$ (in the sense of Definition \ref{our_action}) on ${\rm Bimod}(\cla)$, which agrees with $\Phi$ when restricted on ${\rm Rep}(C^*(S_3))$, is monoidally equivalent to the DQG $C^*(S_3)$. Furthermore, given any DQG $\clq$ with an outer action $\Gamma:{\rm Rep}(\clq)\to {\rm Bimod}(\cla)$, such that $ {\rm Out}(\cla)\subseteq {\rm Im}(\Gamma)$, is either monoidally equivalent to Out$(\cla)$ or monoidally equivalent to $C^*(S_3)$ itself.
\end{abstract}
\section{Introduction}

Quantum groups are well-known symmetry objects in mathematics and physics, employed to study or understand the quantum symmetries of various non-commutative spaces. They mainly originated in the algebraic setting from the pioneering work of Drinfeld and Jimbo (\cite{drinfeld}
,\cite{jimbo}) and in the functional analytical setting in the works of Woronowicz (\cite {woro1},\cite{woro2}) as compact quantum groups, Van-Daele (\cite{dqg}) as discrete quantum groups, locally compact quantum groups in the works of Kustermans and Vaes (\cite{vnqgp}) and as multiplicative unitaries in the works of Baaj and Skandalis (\cite{baaj_skandalis}), just to mention a few of the ever expanding and vast literature. It was only but obvious to try to realise quantum groups as symmetry objects of mathematical structures and already a volume of work has been done in this direction as well starting from Manin(\cite{manin_book}), Wang (\cite{wang1}\cite{wang2}, Banica, Bichon (\cite{ban},\cite{bichon},\cite{bichon_banica}) and many more. Presently, quantum symmetries in the context of discrete and locally compact quantum groups have also garnered attention among mathematicians as is evident form the works of \cite{voigt},\cite{vaes_inf}. The authors have also introduced the notion of discrete quantum groups of automorphisms of (unital) $C^*$ algebras in their work \cite{dg_ss}. 

On the other hand, it is needless to say that category theory has also played an extremely crucial role in the development of the theory of quantum groups as depicted in the work of Woronowicz \cite{woro_tannaka} who proved a quantum analogue of the celebrated Tannaka-Krein Duality theorem. Therefore, in principle, all the properties of a compact quantum group (CQG) and equivalently of the dual discrete quantum group (DQG)  can be formulated from its corepresentation or representation categories and the associated fiber functors respectively. In the same spirit, one can try to understand the coaction of a compact quantum group $\hat{\clq}$ or the dual discrete quantum group $(\clq)$ on a von Neumann algebra, $\cla$ with a finite faithful trace $\tau$, simply by trying to understand some categorical analogue of action of the corresponding Corep($\hat{\clq})$ or, ${\rm Rep}(\clq)$. Indeed, given such a coaction, one can construct a unitary tensor functor $\clf:{\rm Rep}(\clq)\to {\rm Bimod}(\cla)$ and such a functor is called the action of the unitary tensor category ${\rm Rep}(\clq)$ on the von Neumann algebra $\cla$, by the authors of \cite{corey_evans}.
\\
This functor $\clf$ actually lands in the subcategory ${\rm Bimod}_{\mathds{N}}(\cla)$ consisting of $\cla-\cla$ bimodules ${}_{\cla}\clh_{\cla}$, such that ${\rm dim}{}_{\cla}\clh={\rm dim}\clh_{\cla}\in \mathds{N}$(see \cite{krp_insa}). Moreover, $\clf$ preserves dimension, in the sense that, for a representation $\al$ on $\mathds{C}^d$ we have ${\rm dim}{}_{\cla}\clh={\rm dim}\clh_{\cla}=d$, where $\clf(\al)=\clh_{\al}\in {\rm Bimod}_{\mathds{N}}(\cla)$.\\ 
In fact, categorical action in this sense can be interpreted as a generalized outer symmetry for the following reasons. Consider ${\rm Aut}(\cla)$ as a discrete group or rather as the DQG $\clq$ and a functor $\clf:{\rm Rep}(\clq)\to {\rm Bimod}(\cla)$ discussed in the previous paragraph. This sends $\al\in {\rm Rep}(\clq)$ to the bimodule $L^2(\cla,\tau)$ with the trivial right action and the left action given by $b.a:=\al(b)a$ for all $a,b\in \cla$. It is easy to see that $\clf(\al)\cong\clf(\bta)$ if $\al$ and $\bta$ are inner equivalent. In other words, $\clf(\cla)$ only depends on the outer class of $\al$ and $\clf$ descends to a functor from ${\rm Rep}({\rm Out}(\cla)$ to ${\rm Bimod}(\cla)$ which is injective on objects. In view of this, it is natural to think of any unitary tensor functor from ${\rm Rep}(\cls)$ for any DQG $\cls$ to ${\rm Bimod}_{\mathds{N}}(\cla)$ which is also dimension preserving in the sense already discussed, to be a generalization of outer symmetry given by (discrete) quantum group. This is the main philosophy of the present work. We also impose the additional condition of injectivity (at the object level) on the categorical notion as an analogue of faithfulness of usual action or coaction of a quantum group. 
A great deal of work along these lines has been done in the recent times (\cite{corey_evans},\cite{kitamura} and the references therewithin) to name a few. For example, in \cite{corey_evans}, the authors have studied the actions of unitary fusion categories (unitary tensor categories with finitely many isomorphism classes of simple objects) on the non-commutative tori.
\\Thinking along the same lines, one can try to study and possibly classify the action of unitary tensor categories on finite dimensional von Neumann algebras as well. Given a finite dimensional von Neumann algebras $\cla$, the first author in his paper \cite{krp_insa} defined an action of the unitary tensor category ${\rm Rep}(\clq_{\rm aut}(\cla))$ on ${\rm Bimod}(\cla)$ by a tensor functor $\Phi$. But in general, that action may not be injective in the sense of \cite[Definition 4.1 (2)]{kitamura}. In this paper we have tried to understand the quantum symmetries of the von Neumann algebra $\cla:=\C\oplus \C\oplus \M_2(\C)$, via an injective action of the representation category of a DQG. Since $\C\oplus\C\oplus \M_2(\C)$ is the underlying von Neumann algebra of the DQG $C^*(S_3)$, it has an injective action on ${\rm Bimod}(\cla)$ given by the restriction of the functor $\Phi$ and moreover the fusion rules of the ${\rm Rep}(S_3)$ can be employed to show that any DQG $\clq$, containing $C^*(S_3)$ as a quantum subgroup and having an injective action $\tilde{\Phi}$ on ${\rm Bimod}(\cla)$ which is compatible with the action $\Phi$, has to be monoidally equivalent to $C^*(S_3)$ [Theorem \ref{maximal_S3}].
The key step in the proof of our results is the connection between ${\rm Bimod}(\cla)$, for a finite dimensional von Neumann algebra $\cla$, and the matrices of non-negative integer entries as formulated by the first author in \cite{krp_insa}. In this paper, we have dealt with $\C\oplus\C\oplus\M_2(\C)$ but we believe that similar techniques may be employed to obtain the quantum outer symmetry, in the categorical sense, of other finite dimensional von Neumann algebras as well.
\noindent The contents of this paper is briefly outlined below.  
\\
We start by recalling the basic definitions of compact and discrete quantum groups along with their actions on von Neumann algebras and some preliminaries of unitary tensor categories in the first two subsections of section 2. The first part of subsection 2.3 is dedicated to a brief recollection of the generalities of the bimodule category of finite tracial von Neumann algebras. In the second part of subsection 2.3 the map $\Psi$ connecting the bimodules over finite dimensional von Neumann algebras with the square matrices over $\mathds{Z}_{\geq0}$ is stated, the more detailed description of which can be found in \cite{krp_insa}. The last part of section 2 is a brief recollection of the notion of the DQG of automorphisms of a von Neumann algebra, $\cla$, which we denote by $\clq_{\rm aut}(\cla)$. \\
In the third section of our paper we motivate and state the notion of outer action of UTC on bimodule categories in our sense in Definition \ref{our_action}. 
\\The fourth section contains the main results of this paper. We start by fixing our finite dimensional von Neumann algebra $\cla$ as $\C\oplus\C\oplus\M_2(\C)$, which by the virtue of being the underlying von Neumann algebra of the discrete quantum group $C^*(S_3)$ has a coaction of $C^*(S_3)$ on itself. Therefore, $C^*(S_3)$ becomes a quantum subgroup of $\clq_{\rm aut}(\cla)$ and hence we can restrict the functor $\Phi$ on ${\rm Rep}(C^*(S_3)$. Then we have the following 
\bthm [Theorem \ref{maximal_S3}]
Let $\tilde{\Phi}:{\rm Rep}(\clq)\to {\rm Bimod}(\cla)$ be an outer action such that $\tilde{\Phi}\circ\tilde{\pi}=\Phi.$ Then ${\rm Rep}(\clq)$ is equivalent to ${\rm Rep}(C^*(S_3))$, i.e. the two DQGs $\clq$ and $C^*(S_3)$ are monoidally equivalent.
\ethm 
\bcrlre [Corollary \ref{contains_Out}]
    Let $\clq$ be any DQG with an outer action $\Gamma:{\rm Rep}(\clq)\to {\rm Bimod}_{\mathds{N}}(\cla)$ such that ${\rm Im}(\tilde{\Phi})$ contains ${\rm Out}(\cla)$. Then $\clq$ is monoidally equivalent to either ${\rm Out}(\cla)$ or $C^*(S_3)$. \ecrlre
  Finally, invoking Corollary 1.4 of \cite{etingof_gelaki}, we conclude the paper with 
  \bthm [Theorem \ref{grp_algebra}]      Let $\clq=C^*(H)$ for a finite group $H$, such that ${\rm Rep}(\clq)$ satisfies Theorem \ref{maximal_S3} and Corollary \ref{contains_Out}. Then, $H$ is isomorphic to $S_3$ and hence $\clq\cong C^*(S_3)$.
    \ethm.
    \section{Preliminaries} 
\subsection{Compact and discrete quantum groups (von Neumann algebraic sense) and their coaction on von Neumann algebras}
We shall work in the framework of von Neumann algebraic compact and discrete quantum groups. We shall use the symbol $\otimes$ to denote both the algebraic and von Neumann algebraic tensor products unless there is any scope of confusion. 
\bdfn\label{vna_cqg}
A von Neumann algebraic compact quantum group (CQG) is a von Neumann algbera $\hat{\clq}$ equipped with a unital, normal $\ast$ homomorphism $\hat{\Delta}:\hat{\clq}\to \hat{\clq}\otimes \hat{\clq}$ along with a faithful state $h$ which is bi-invariant in the sense that $({\rm id}\otimes h)\circ\hat{\Delta}(a)=h(a)1=(h\otimes {\rm id})\circ\hat{\Delta}(a)$ for all $a\in \hat{\clq}$.\edfn
\bdfn\label{vna_dqg}A von Neumann algebraic discrete quantum group $\clq$ is given by the von Neumann algebraic closure of a direct sum of finite dimensional matrix algebras , i.e. $\clq=\prod_{i\in I}\clb(\clh_i)$, where $I$ is an index set, equipped with a unital $\ast$ normal homomorphism $\Delta:\clq\to \clq\otimes \clq$ and a pair of faithful semi finite weights $h_L$ and $h_R$, having the finite dimensional matrix algebras $\clb(\clh_i)$ in their domains and satisfying suitable  left and right invariance conditions respectively.\edfn
\noindent Compact and discrete quantum groups, both belong to the category of locally compact quantum groups and have a dual pairing between them. For more details one can refer to \cite{vnqgp}. 
\\Analogous to the action of groups on von Neumann algebras, we have the notion of coaction of quantum groups on von Neumann algebras as follows
\bdfn\label{coaction_vNa}
\begin{enumerate}
\item 
Given a von Neumann algebra $\cla$, a coaction of a DQG $(\clq,\Delta)$ on $\cla$ is given by a unital, coassocitive, normal and injective $\ast$ homomorphism $\al:\cla\to \cla\otimes \clq$ such that the linear span of $\al(\cla)(1\otimes \clq)$ is ultra-weakly dense in $\cla\otimes \clq$. 
\item A von Neumann algebraic coaction $\al:\cla\to \cla\otimes \clq$ is called faithful if the linear span of $(\omega\otimes {\rm id})\circ\al(a):\omega\in \cla_*,a\in \cla\}$ is ultra-weakly dense in $\clq$.\end{enumerate}
\edfn
\subsection{Unitary tensor categories}
 We state a few basic examples of unitary tensor categories (UTC for short) and recall a few standard definitions in this subsection.
 The two standard examples of UTCs which are central to this paper are the corepresentation category of a CQG and the representation category of a DQG defined as follows.
\bdfn\label{Corep}
Given a compact quantum group $\hat{\clq}$, we denote the UTC of finite dimensional unitary corepresentations of $\hat{\clq}$ by ${\rm Corep}(\hat{\clq})$. The morphisms in this category are given by the intertwiners between corepresentations.  \edfn
\bdfn\label{Rep}Given a DQG $\clq$, we have the unitary tensor category of finite dimensional $\ast$ representations of the DQG $\clq$, namely, ${\rm Rep}(\clq)$. \edfn
\brmrk\label{forget_fiber}
${\rm Corep}(\hat{\clq})$ is naturally equipped with a fiber functor which sends every object to the corresponding finite dimensional complex Hilbert space associated with it. Similarly, ${\rm Rep}(\clq)$ is equipped with the natural forgetful fiber functor.\ermrk 
Finally, we have 
\bthm(Woronowicz's Tannaka-Krein Duality)\cite{woro_tannaka}\cite [Theorem 2.3.2]{sergey_book}
Let $\mathscr{C}$ be a unitary tensor category with conjugates, $F:\mathscr{C}\to \text{Hilb}_{f}$ be a unitary fiber functor. Then there exists a compact quantum group $\hat{\mathcal{Q}}$ and a unitary monoidal equivalence $E:\mathscr{C}\to \text{Corep}(\hat{\mathcal{Q}})$ such that $F$ is naturally monoidally isomorphic to the composition of the canonical fiber functor $\text{Corep}(\hat{\mathcal{Q}})\to \text{Hilb}_{f}$ with $E$. Further more, the Hopf $*$- algebra $(\hat{\mathcal{Q}}_0,\Delta)$ is uniquely determined upto isomorphism. In case, $\mathscr{C}$ is of Kac-type we get $\hat{Q}$ to be of Kac-type as well.
\ethm
This is an extremely crucial tool to construct compact and discrete quantum groups from UTCs equipped with fiber functors. 
\\It is worth mentioning at this point that not all $C^*$ tensor categories admit fiber functors \cite[Example 2.3.7]{sergey_book} and a $C^*$-tensor category can also have fiber functors producing nonisomorphic quantum groups. This leads us to the following 
\bdfn\label{monoidal_equiv}
Compact quantum groups $\hat{\clq}_1$ and $\hat{\clq}_2$ are said to be monoidally equivalent if the categories ${\rm Corep}(\hat{\clq}_1)$ and ${\rm Corep}(\hat{\clq}_2)$ are unitarily monoidally equivalent.\edfn 
Two dual DQGs $\clq_1$ and $\clq_2$ are called monoidally equivalent if their dual CQGs are monoidally equivalent.  
\brmrk\label{example_mon}
A genuine compact group maybe monoidally equivalent to a genuine compact quantum group as is portrayed in Example 2.3.9 of \cite{sergey_book}.\ermrk

We refer to \cite{sergey_book} and \cite{etingof_fusion}, for detailed discussion on UTCs and Tannaka-Krein reconstruction.
\subsection{Bimodule categories of tracial von Neumann algebras}
\subsubsection{Generalities of bimodule categories of von Neumann algebras}
Let $\cla$ be a tracial von Neumann algebra with a finite faithful tracial state $\tau$. Then we have the notion of a $\cla-\cla$ bimodule ${}_{\cla}\clh_{\cla}$, which is given by a Hilbert space $\clh$ and left and right actions of $\cla$ on $\clh$. We have the following notion of Connes' tensor product of two bimodules over $\cla$. 
\bdfn\label{Popa_Ananth_Connes}[Definition 13.2.4 \cite{popa_ananth}]
Let ${}_{\cla}\clh_{\cla}$ and ${}_{\cla}\clk_{\cla}$ be $(\cla,\tau)$ bimodules. The Connes tensor product of ${}_{\cla}\clh_{\cla}$ and ${}_{\cla}\clk_{\cla}$, denoted by $\clh\otimes_{\cla}\clk$, is given by the separation and completion of the vector space $\clh^0\otimes _{\rm alg}\clk$ with respect to the sesquilinear form $\langle\xi_1\otimes_{\rm alg} \xi_2,\eta_1\otimes_{\rm alg}\eta_2\rangle:=\langle\eta_1,\langle\xi_1,\xi_2\rangle_{\cla}\eta_2\rangle_{\clk}$ defined on $\clh^0\otimes _{\rm alg}\clk$, where $\clh^0$ denotes the set of all left-bounded vectors in the right module $\clh_{\cla}$ as in \cite{popa_ananth}.\edfn
For equivalent definitions of Connes tensor product of bimodules for which we refer to \cite{popa_ananth}. \\
Give such a bimodule over a tracial von Neumann algebra $(\cla,\tau)$, we have the notion of a left and right dimensions as is given in the following
\bdfn\label{dimension_of_module}
Let $\clh_{\cla}$ is a right $\cla$ module over a tracial von Neumann algebra $(\cla,\tau)$. The scalar $({\rm Tr}\otimes \tau)(p)$ where $p$ is any projection in $(\clb(l^2(\mathds{N})\otimes \cla)$ such that $\clh_{\cla}\cong p(l^2(\mathds{N})\otimes L^2(\cla,\tau))$, denoted by ${\rm dim}(\clh_{\cla})$, is called the right dimension/rank of $\clh_{\cla}$. Similarly, we have the notion of a left dimension/rank of a left module ${}_{\cla}\clh$.\edfn 
We therefore, have the $C^*$ multi-tensor category ${\rm Bimod}(\cla)$ whose objects are Hilbert $\cla-\cla$ bimodules on $\cla$, ${}_{\cla}\clh_{\cla}$, which are bi-finite i.e. ${\rm dim}(\clh_{\cla})<\infty$ and ${\rm dim}({}_{\cla}\clh)<\infty$. The morphisms are given by $\cla-\cla$ bimodule maps and the monoidal structure given by the Connes tensor product of bimodules, with ${}_{\cla}L^2(\cla,\tau)_{\cla}$ being the unit object. Furthermore, given any bi-finite index bimodule ${}_{\cla}\clh_{\cla}$, we have the contragradient bimodule ${}_{\cla}\overline{\clh}_{\cla}$ with the left action given by $a.\overline{h}:=\overline{h.a^*}$ and the right action given by $\overline{h}.a:=\overline{a^*.h}$. ${}_{\cla}\overline{\clh}_{\cla}$ also becomes a dual object to ${}_{\cla}\clh_{\cla}$ in the bimodule category, thereby making ${\rm Bimod}(\cla)$ into a rigid $C^*$ multi tensor category. [ In particular, if $(\cla,\tau)$ is a tracial factor, then ${\rm Bimod}(\cla)$ becomes a rigid $C^*$ tensor category.]
\\As mentioned in the introduction, there is an important sub-category of ${\rm Bimod}(\cla)$ which plays a crucial role in our work. 
\bdfn\label{finite_rank_subbi}
 ${\rm Bimod}_{\mathds{N}}(\cla)$ denotes the full-subcategory of $\cla-\cla$ bimodules, ${}_{\cla}\clh_{\cla}$ such that ${\rm dim}{}_{\cla}\clh={\rm dim}\clh_{\cla}\in \mathds{N}$.We define this integer as the rank of the bimodule ${}_{\cla}\clh_{\cla}\in {\rm Bimod}_{\mathds{N}}(\cla)$. \edfn

 We end this subsection with the following
 \bthm\label{Out(A)insideGrp(A)}
Given any tracial von Neumann algebra $(\cla,\tau)$, the outer automorphism group of $\cla$, ${\rm Out}(\cla)$ embeds inside ${\rm Grp}(\cla)$, where ${\rm Grp}(\cla)$ consists of the irreducible elements ${}_{\cla}\clh_{\cla}$ in ${\rm Bimod}(\cla)$ such that ${}_{\cla}\clh\otimes _{\cla}\overline{\clh}_{\cla}\cong {}_{\cla}L^2(\cla,\tau)_{\cla}$.
\ethm
We refer the reader to \cite{popa_ananth}, \cite{bimod_1}, \cite{bimod_2}, \cite{bimod_3} and \cite{BDH} for a detailed discussion of the notions above.
\subsubsection{Bimodule category of finite dimensional matrix algebras and their connection with matrices with non-negative integer entries}
\noindent Let $\cla=\M_{d_1}(\C)\oplus\M_{d_2}(\C)\oplus...\oplus\M_{d_k}(\C)$.
    Let us write $\cla_i=\M_{d_i}(\C)$ with $p_i$ being the central projection corresponding to $\cla_i$,i.e. $p_i\cla=\cla_i$.
    \\Let $a_l$ and $a_r$ denote the left and right actions of $a\in \cla$ on $L^2(\cla,Tr)$.\\
    As any $\ast$-homomorphism between full matrix algebras is determined by the multiplicity matrix upto unitary equivalence, the bimodule category Bimod$(\cla)$ is a multi-tensor category with the irreducibles $E_{ij}=\C^{d_i}\otimes \overline{\C^{d_j}}$ with the left and right actions of $\cla$ given by $\M_{d_i}(\C)$ on $\C^{d_i}$ and $\M_{d_j}(\C)$ on $\overline{\C^{d_j}}$ respectively. Thus any bimodule upto isomorphism is determined by the multiplicities $m_{ij}$ of $E_{ij}$.
    Let us recall from \cite{krp_insa} that if a bimodule is given by a Hilbert space $\clh$, with the left action of $\cla$ given by the $*$ homomorphism $\al:=\pi_l:\cla\to \clb(\clh)$ and the right action given by $\pi_r$, then $m_{ij}=\frac{Tr(\pi_l(p_i)\pi_r(p_j))}{d_id_j}$ where $Tr$ denotes the usual matrix trace.
    We also refer to \cite{krp_insa} to recall that if $\overline{X}$ is the the dual object of $X$, then the corresponding multiplicities of $E_{ij}$ in $\overline{X}$ given by $\overline{m_{ij}}$ can be calculated as $\overline{m_{ij}}=m_{ji}$.
    \\The preceeding discussion can be written in the form of the following theorem [Thm 3.6 , \cite{krp_insa}]
    \bthm\label{thm_for_Psi}
    Let $\cla=\M_{d_1}(\C)\oplus\M_{d_2}(\C)...\oplus\M_{d_k}(\C)$ as considered in the beginning of this subsection. There is  bijective map $\Psi$ from Obj(Bimod($\cla))$ onto $M_{k\times k}(\mathds{Z})$, the set of $k\times k$ matrices with non negative integer entries such that $\Psi(X)=(m^X_{ij})$ such that $m_{ij}^X$ is the multiplicity of $E_{ij}$ in $X$. This is a $Z_{\geq 0} $ ring isomorphism with $\Psi(\overline{X})=\Psi(X)'$ where $'$ denotes the usual matrix transpose. 
    \ethm
     Now we restrict ourselves to the sub category $\rm Bimod_{\mathds{N}}(\cla)$ and see that the objects in this subcategory map under $\Psi$ to the matrices a particular form as is described in the following :
     \bcrlre\label{image_Psi}[Corollary 3.7 \cite{krp_insa}]
     Under the isomorphism $\Psi$ Obj($\rm Bimod_{\mathds{N}}(\cla)$ gets mapped onto the subset $\mathscr{F}$ of $\M_{k\times k}(\mathds{Z})$ consisting of $k\times k$ matrices $A$ with non negative integer entries such that both $A$ and $A'$ have the vector $d=(d_1,d_2,...,d_k)$ as an eigen vector corresponding to a common non-negative integral eigen value i.e. $\mathscr{F}=\bigcup_{n>0}\{A\in \M_{k\times k}(\mathds{Z}):Ad=nd\text{ and }A'd=nd\}.$
     \ecrlre
    \brmrk\label{definition of rank}
     We call the rank of the matrix $A$, satisfying $Ad=nd$ and $A'd=nd$ to be $n$. Observe that the a bimodule of rank $n$ in ${\rm Bimod}(\cla)$ gets mapped to a matrix of rank $n$ in $\clf$. Hence if we start with a $n$ dimensional representation $(\al,\clh)$ in ${\rm Rep}(\clq_{\rm aut}(\cla))$, then the rank of $\Psi(\Phi(\al,\clh_{\al})=n$.\ermrk
     Given any two objects $A$ and $B\in \mathscr{F}$, we call $A\leq B$ if $A_{ij}\leq B_{ij}\text{for all }i,j\in \{1,2,...,k\}.$ We call an element $A$ to be irreducible if $A\geq B$ for any $B\in \mathscr{F}$ implies $A=B$. \\
      We end this subsection with following 
 \blmma\label{general_lemma}
 Given any finite dimensional von Neumann algebra $\cla$ the set $Im(\Psi\circ\Phi)\subseteq\mathscr{F}_{\cla}$ is closed under sum, product and transpose.
 \elmma
 \subsection{Discrete quantum group of automorphisms of a von Neumann algebra}
Given a von Neumann algebra $\cla$, consider the category $\widetilde{\cld_{\cla}}$ (as introduced by the first author in \cite{krp_insa}) whose objects are given by $X=(\alpha, \clh)$, where $\clh$ is a finite dimensional complex Hilbert space, $\alpha : \cla \raro \cla \ot \clb(\clh)$ is a  unital $\ast$-homomorphism. Morphisms  ${\rm Mor}(X,Y)$ for  objects $X=(\alpha, \clh)$ and $Y=(\beta,  \clk)$ are given by  linear maps $T : \clh \raro \clk$  
      such that $$(1 \ot T) \alpha(a)=\beta(a) (1 \ot T) ~ \forall a \in \cla.$$   Tensor product and direct sum operations are as usual, given by
       $$ (\alpha, \clh) \oplus (\beta, \clk):=(\alpha \oplus \beta, \clh \oplus \clk) \text { and } (\alpha, \clh) \otimes (\beta, \clk):=((\alpha \ot {\rm id}) \circ \beta, \clh \ot \clk) \text{ respectively} .$$ The unit object is given by $\mathds{1}:=({\rm id}_\cla, \IC)$. It is easy to see that $\mathds{1} \ot X \cong X \cong X \ot \mathds{1}$ for any object $X\in \widetilde{\cld_{\cla}}$. For any projection $P\in {\rm Mor}(X,X)$ we can define the sub-object $X_P=(\al_P,\clh_P)$ where $\al_P(a)=(\rm id_{\cla}\otimes P)\al(a)(\rm id_{\cla}\otimes \iota_{P\clh})$ where $\iota_{P\clh}:P\clh\to \clh$ is the inclusion map.  Therefore, $\widetilde{\cld_{\cla}}$ is a $C^*$ tensor category or unitary tensor category (UTC) which is closed under taking finite direct sums, tensor products and sub-objects.
    \\We call an object $(\al,\clh_{\al})\in \widetilde{\cld_{\cla}}$ to be dualizable if the following conditions hold:\\
     there are $s \in \overline{\clh} \ot \clh$, $t \in \clh \ot \overline{\clh}$ (where $\overline{\clh}$ denotes the complex conjugate Hilbert space of $\clh$), $\hat{\alpha} : \cla \raro \cla \ot \clb(\overline{\clh})$  unital $\ast$-homomorphism satisfying
     $$( (\hat{\alpha}  \ot {\rm id}) \circ \alpha(a)) (1 \ot s)=a \ot s,~~~((\alpha \ot {\rm id}) \circ \hat{\alpha}(a))(1 \ot t)=a \ot t~~~\forall a \in \cla,$$
     $$ (R^*_s \ot 1_{\overline{\clh}})(1_{\overline{\clh}} \ot R_t)=1_{\overline{\clh}},~~~~(R^*_t \ot 1_\clh)(1_\clh \ot R_s)=1_\clh,$$
      where for a vector $v $ in some finite dimensional vector space $\clk$ we denote by $R_v$ the map from $\IC$ to $\clk$ given by $z \mapsto zv $. 
\\We define $\cld_\cla$ to be the subcategory of all dualizable objects in $\widetilde{\cld_{\cla}}$. We have the following theorem by Longo and Roberts 
\bthm\label{longo_roberts}[Theorem 2.4, \cite{longo_roberts}]
Let $\widetilde{\clc}$ be a $C^*$ tensor category and let $\clc$ be the subcategory consisting of all dualizable objects. Then $\clc$ is closed under conjugates and tensor products. If $\widetilde{\clc}$ has finite direct sums and sub-objects, then $\clc$ is also closed under finite direct sums and sub-objects.\ethm Therefore, $\cld_{\cla}$ being a  rigid full subcategory of $\widetilde{\cld_{\cla}}$ is also closed under tensor product, finite direct sum and sub-objects.    
This gives us the following
\bthm\label{tannaka_krein}[Theorem 2.2\cite{krp_insa}]
$\cld_{\cla}$ is a UTC with a canonical fiber functor $F:\cld_{\cla}\to {\rm Hilb}_f$ given by $F((\al,\clh))=\clh$ and $F(T)=T$ for every morphism $T$.\ethm
Therefore, we can invoke Tannaka-Krein duality theorem to obtain a DQG, say, $\clq_{\rm aut}(\cla)$ such that ${\rm Rep}(\clq_{\rm aut}(\cla)$ is equivalent to $\cld_{\cla}$. Moreover, we have the following universal property of $\clq_{\rm aut}(\cla)$
\bthm\label{universaility}[Theorem 2.3,\cite{krp_insa}]
    (i) There is a faithful coaction $\Theta:\cla\to \cla\otimes \clq_{\rm aut}(\cla)$.\\
    (ii) For any DQG $\cls$, with a coaction $\delta:\cla\to \cla\otimes \cla$, we have a normal unital $\ast$ homomorphism $\hat{\delta}:\clq_{\rm aut}(\cla)\to \cls$, such that $({\rm id}\otimes \hat{\delta})\circ\Theta=\delta.$\\
    (iii) If the coaction $\delta$ is faithful, the the map $\hat{\delta}$ is surjective.

\ethm 
\begin{proof}
        (i) The map $\Theta$ has been defined explicitly in \cite{dg_ss} and the fact that it is a faithful coaction follows easily from the proofs of Theorem 3.2 and Theorem 3.4 of \cite{dg_ss}.\\
     (ii)The proof follows from a verbatim adaptation of the proof of Theorem 3.8 (i) of \cite{dg_ss}.\\
     (iii)The proof follows from the proof of Theorem 3.8 (ii) of \cite{dg_ss}. 
\end{proof}

\section{ Categorical action of UTC on ${\rm Bimod}(\cla)$}
We recall the following definition given by \cite{corey_evans} of categorical action in
\bdfn\label{action_tensor}(\cite{corey_evans }, Definition 2.1)
An action of a unitary tensor category $\clc$  on a von Neumann algebra is a $C^*$ tensor functor $\clf:\clc\to {\rm Bimod}(\cla)$.\edfn
As mentioned in the introduction, given the coaction of a DQG $\clq$ on a von Neumann algebra $\cla$, we can get a categorical action in the above sense with certain additional properties. 
Let $\delta:\cla\to \cla\otimes \clq$ be an action of a DQG $\clq:=\prod_{i\in I}\clb(\clh_i)$ on the von Neumann algebra $\cla$. Then, we can write $\delta=\prod_{i\in I}\delta_i$, where $\delta_i:\cla\to \cla\otimes \clb(\clh_i)$, is a unital $\ast$-homomorphism. Then, we can define the functor \begin{equation}\label{eqn}\clf\equiv\Phi_{\delta}:{\rm Rep}(\clq)\ni(\rho,\clh_{\rho})\to L^2(\cla)\otimes \clh_{rho}\in {\rm Bimod}(\cla),\end{equation} with the trivial right action on $L^2(\cla,\tau)\otimes \clh_{\rho}$ and the left action given by $b.(a\otimes \xi):=(({\rm id}\otimes \rho)\circ\delta)(b)(a\otimes \xi)$ for all $a,b\in \cla $ and $\xi\in \clh_{\rho}$. 
\brmrk\label{Sir's functor}
\begin{enumerate}
    \item
One can take $\clq=\clq_{\rm aut}(\cla)$ with the canonical coaction on $(\cla,\tau)$ as in subsection 2.4. Then the map \begin{equation}\label{sir_func}\Phi:\cld_{\cla}\ni (\al,\clh_{\al})\to L^2(\cla,\tau)\otimes \clh_{\al}\in {\rm Bimod}(\cla), \end{equation} with the trivial right action and the left action given by $\al$,defines a tensor functor such that $\Phi$ is identity on the morphism spaces\cite{krp_insa}.
\item 
 Given any element $\al\in {\rm Aut}(\cla)$, $\Phi(\al,\mathds{C})\in {\rm Grp}(\cla).$\end{enumerate}
 \ermrk
 \bdfn\label{dim_pre}
 We call an action $\clf$ of ${\rm Rep}(\cls)$ for some DQG $\cls$ on ${\rm Bimod}(\cla)$, to be dimension preserving if the range of $\clf$ is in ${\rm Bimod}_{\mathds{N}}(\cla)$ and $${\rm dim}{}_{\cla}\clf(\rho,\clh_{\rho})={\rm dim}{}_{\cla}\clf(\rho,\clh_{\rho})={\rm dim}(\rho,\clh_{\rho}).$$
 \edfn
 As discussed in the introduction, in case of the DQG associated with ${\rm Aut}(\cla)$, with the canonical coaction, say $\delta$, the functor ${\Phi_{\delta}}$ descends to a functor on ${\rm Rep}({\rm Out}(\cla))$, which is injective on objects. This motivates us to define a categorical version of outer quantum symmetry in the following 
 \bdfn\label{our_action}
 A categorical outer symmetry of $(\cla,\tau)$ given by a DQG $\cls$ is a unitary tensor functor $\clf:{\rm Rep}(\cls)\to {\rm bimod}_{\mathds{N}}(\cla)$ which is dimension preserving and injective on objects. We will often call such an action to be an outer action of ${\rm Rep}(\cla)$on $(\cla,\tau)$. \edfn
 \brmrk\label{contrast_to_kitamura}
 We should mention at this point that an outer action in our sense (Definition \ref{our_action}) is not necessarily full and hence is not objectwise outer in the sense of Definition 4.1 (3) in \cite{kitamura}.\ermrk
 The aim of this paper is to explore such categorical outer symmetry for the special case of $\cla=\C\oplus \C\oplus \M_2(\C)$, which is also the underlying von Neumann algebra of the DQG $C^*(S_3)$.

\section{Main Result}
We fix $\cla=\mathds{C}\oplus \mathds{C}\oplus \M_2(\mathds{C})$ as the main object of study in this section. We refer to Example 3 in subsection 3.3 of \cite{krp_insa}, which explains that $C^*(S_3)$ can be identified as a quantum subgroup of $\clq_{\rm aut}(\cla)$.

\subsection{Outer action of ${\rm Rep}(\clq)$ on $\cla$ such that ${\rm Rep}(C^*(S_3))\subseteq {\rm Rep}(\clq)$}
\noindent 
Let $\clq$ be a discrete quantum group which contains $C^*(S_3)$ has a quantum subgroup. Then there is a surjective discrete quantum group homomorphism $\pi:\clq\to C^*(S_3)$ which induces the functor $\tilde{\pi}:{\rm Rep}(C^*(S_3))\to {\rm Rep}(\clq)$ such that for any $\al\in {\rm Rep}(C^*(S_3))$, $\tilde{\pi}(\al)(q):=\al\circ\pi(q)$ for all $q\in \clq$. Recall from \cite{krp_insa} the bimodules  $E_1=_{\al_1}L^2(\cla,\tau)\otimes {\clh_{\al_1}}, E_2=_{\al_2}L^2(\cla,\tau)\otimes {\clh_{\al_2}}\text{ and }E_3=_{\al_3}L^2(\cla,\tau)\otimes{\clh_{\al_3}}$ corresponding to the objects $\al_1,\al_2$ and $\al_3\in {\rm Rep}(C^*(S_3)$, where $\cla$ acts from the right via the trivial action and on the left by the morphisms $\al_i$. Since the Obj(${\rm Rep}(C^*(S_3))\cong \mathds{Z}_{\geq0}\{\al_1,\al_2,\al_3\}$, hence we can restrict the functor $\Phi$, introduced in Remark \ref{Sir's functor}, to ${\rm Rep}(C^*(S_3))$ which is a full subcategory of $\cld_{\cla}^{\kac,\tau}$ to obtain  $\Phi:{\rm Rep}(C^*(S_3))\to {\rm Bimod}(\cla)$ which takes $\al_i$ to the corresponding bimodule.  Then we have the following 
\bthm\label{maximal_S3}
Let $\tilde{\Phi}:{\rm Rep}(\clq)\to {\rm Bimod}(\cla)$ be an outer action (in the sense of Definition \ref{our_action}) such that 
the following diagram 
\[\begin{tikzcd}
	{{\rm Rep}(\mathcal{Q})} & {{\rm Bimod}(\mathcal{A})} \\
	{{\rm Rep}(C^*(S_3))}
	\arrow["{\tilde{\Phi}}", from=1-1, to=1-2]
	\arrow["{\tilde{\pi}}", from=2-1, to=1-1]
	\arrow["{\Phi}", from=2-1, to=1-2]
\end{tikzcd}\]commutes, i.e. $\tilde{\Phi}\circ \tilde{\pi}=\Phi.$ Then ${\rm Rep}(\clq)$ is equivalent to ${\rm Rep}(C^*(S_3))$, i.e. the two DQGs $\clq$ and $C^*(S_3)$ are monoidally equivalent.
\ethm
We prove this theorem through a series of computations.
According to the notations used in subsection 3.2, primarily in Theorem \ref{thm_for_Psi} and Corollary \ref{image_Psi}, for $\cla=\C\oplus\C\oplus \M_2(\C)$, we have $k=3,d_1=1=d_2\text{ and }d_3=2$ i.e. $d=(1,1,2)$ 
\\One can easily show that $A_1=\begin{pmatrix}1&0&0\\0&1&1\\0&0&1\end{pmatrix}~A_2=\begin{pmatrix}0&1&0\\1&0&0\\0&0&1\end{pmatrix}~A_3=\begin{pmatrix}0&0&1\\0&0&1\\1&1&1\end{pmatrix}$ satisfy the relations $$Ad=nd\text{ and }A'd=nd\text{ for }n=1,1\text{ and }2\text{ respectively}.$$ Furthermore, from \cite{krp_insa}(we know that
$\Psi\circ{\tilde{\Phi}}\circ\tilde{\pi}(\alpha_i))=A_i\text{ for }i=1,2,3.$ Hence, $A_1,A_2$ and $A_3\in {Im}(\Psi\circ\tilde{\Phi})\subseteq\mathscr{F}_{\cla}.$  
\\Our main aim now is to find the generating set of the image of $\Psi\circ {\tilde{\Phi}}\subseteq \mathscr{F}$. 
\\In that pursuit we first try to follow the following algorithm:
\\Let $A\in Im(\Psi\circ \tilde{\Phi})$ and $m_1$ be the largest non negative integer such that $m_1A_1$ is a sub-object of $A$ i.e. $m_1A_1\leq A$. Similarly let $m_2$ be the largest non negative integer such that $m_2A_2\leq A-m_1A_1$ and $m_3$ be the largest non negative integer such that $m_3A_3\leq A-m_1A_1-m_2A_2.$ 
\\If, for every such $A$ we have $A=m_1A_1+m_2A_2+m_3A_3$, then we are done. On the other hand, let there be some A $\in Im(\Psi\circ\tilde{\Phi})$ such that $A=m_1A_1+m_2A_2+m_3A_3+X$ for some non zero $X\in Im(\Psi\circ\tilde{\Phi}).$
\\In order to understand the structure of such a $X\in Im(\Psi\circ\tilde{\Phi})$, we start with the following

\blmma For any $A\in  Im(\Psi\circ\tilde{\Phi}),$ for $\cla=\mathds{C}\oplus\mathds{C}\oplus\M_2(\mathds{C}),$ $A_{33}\neq 0.$\elmma
\begin{proof}
    Let us assume $a_{33}=0$. Since $Ad=nd\text{ and }A'd=nd$, hence we have,
    \begin{equation*}
      \begin{pmatrix}
        a_{11}&a_{12}&a_{13}\\a_{21}&a_{22}&a_{23}\\a_{31}&a_{32}&0\end{pmatrix}\begin{pmatrix}1\\1\\2\end{pmatrix}=\begin{pmatrix}n\\n\\2n
    \end{pmatrix} \text{ and }\begin{pmatrix}a_{11}&a_{21}&a_{31}\\a_{12}&a_{22}&a_{32}\\a_{13}&a_{23}&0\end{pmatrix}\begin{pmatrix}1\\1\\2\end{pmatrix}=\begin{pmatrix}n\\n\\2n\end{pmatrix}  
    \end{equation*}
  which implies \begin{equation}\label{eqn2}
      a_{31}+a_{32}=2n,
  \end{equation}
  \begin{equation}\label{eqn3}
      a_{11}+a_{21}+2a_{31}=n
  \end{equation} and 
  \begin{equation}\label{eqn4}
      a_{12}+a_{22}+2a_{32}=n.
  \end{equation} Multiplying $\eqref{eqn2}$ by 2 and putting it in $\eqref{eqn3}+\eqref{eqn4},$  we get $$a_{11}+a_{21}+4n+a_{21}+a_{22}=2n$$ which is not possible given all the $a_{ij}\geq0.$
\end{proof} 
Hence, given $X\in \mathscr{F}_{\cla}$, we have $x_{33}\neq 0$. Moreover, by virtue of the algorithm which we have used to obtain $X$, we can say that $X$ cannot dominate $A_1,A_2\text{ and }A_3$ which further restricts the possible structure of $X$ within the following cases: 
\begin{enumerate}
    \item $X \ngeq A_1$ which means $x_{11}=0~ \text{or}~ x_{22}=0$,
    \item $X\ngeq A_2$ which means $x_{12}=0~\text{or} ~x_{21}=0$ and
\item $X\ngeq A_3$ which means $x_{13}=0~or~x_{23}=~\text{or}
~ x_{32}=0.$
\end{enumerate}
 This along with the fact that $Xd=nd\text{ and }X'd=nd$ gives the following possible choices of $X$:\\\\
 
 \noindent $
     (1) \begin{pmatrix}
         0&0&x_{13}\\x_{12}&x_{22}&0\\x_{31}&x_{32}&x_{33}
     \end{pmatrix}~~~
     (2) \begin{pmatrix}
         0&0&x_{13}\\x_{21}&x_{22}&x_{23}\\0&x_{32}&x_{33}
     \end{pmatrix}~~~
     (3) \begin{pmatrix}
         0&0&x_{13}\\x_{21}&x_{22}&x_{23}\\x_{31}&0&x_{33}
         \end{pmatrix}~~~
     (4) \begin{pmatrix}
        x_{11}&x_{12}&x_{13}\\
        0&0&x_{23}\\ 0&x_{32}&x_{33}
    \end{pmatrix}~~~
    (5)\begin{pmatrix}
        x_{11}&x_{12}&x_{13}\\0&0&x_{23}\\x_{31}&0&x_{33}
    \end{pmatrix}\\\\
    (6)\begin{pmatrix}
        x_{11}&x_{12}&0\\0&0&x_{23}\\x_{31}&x_{32}&x_{33}
    \end{pmatrix}~~~
    (7)\begin{pmatrix}
        0&x_{12}&0\\0&x_{22}&x_{23}\\x_{31}&x_{32}&x_{33}
    \end{pmatrix}~~~
    (8)\begin{pmatrix}
        0&x_{12}&x_{13}\\0&x_{22}&0\\x_{31}&x_{32}&x_{33}
    \end{pmatrix}~~~
    (9)\begin{pmatrix}
        0&x_{12}&x_{13}\\0&x_{22}&x_{23}\\0&x_{32}&x_{33}
    \end{pmatrix}
    (10)\begin{pmatrix}
        x_{11}&0&0\\x_{21}&0&x_{23}\\x_{31}&x_{32}&x_{33}
    \end{pmatrix}\\\\
    (11)\begin{pmatrix}
        x_{11}&0&x_{13}\\x_{21}&0&0\\x_{31}&x_{32}&x_{33}
    \end{pmatrix}~~
    (12)\begin{pmatrix}
        x_{11}&0&x_{13}\\x_{21}&0&x_{23}\\0&x_{32}&x_{33}
    \end{pmatrix}$ \\\\
    
   \noindent  Again, since each matrix $X$ satisfies the relations $Xd=nd\text{ and }X'd=nd$ hence we can find a suitable choices of the non negative integers $x_{ij}$ for each of the 12 cases above:
    \begin{enumerate}
        \item For any $X$ of the form (1) we have 2X= $(2k-l)\begin{pmatrix}
            0&0&2\\0&4&0\\2&0&3
        \end{pmatrix}+l\begin{pmatrix}
            0&0&2\\4&0&0\\0&2&3
        \end{pmatrix}$ for $k,l\in \mathds{Z}_{\geq0}$ such that $2k-l\geq0$;
        \item Any $X$ of the form (2) can be written as $k\begin{pmatrix}
            0&0&2\\4&0&0\\0&2&3
        \end{pmatrix} $ for any $k\in \mathds{Z}_{\geq 0} $;
        \item Any $X$ of the form (3) can be written as $k\bpmat 0&0&2\\0&4&0\\2&0&3\epmat$ for any $k\in \mathds{Z}_{\geq0}$;
        \item Any $X$ of the form (4) can be written as $k\bpmat4&0&0\\0&0&2\\0&2&3\epmat$ for any $k\in \mathds{Z}_{\geq0};$
        \item Any $X$ of the form (5) can be written as $k\bpmat0&4&0\\0&0&2\\2&0&3\epmat$ for any $k\in \mathds{Z}_{\geq0}$;
        \item For any matrix $X$ of the form (6) we have $2X=(2k-l)\bpmat0&4&0\\0&0&2\\2&0&3\epmat+l\bpmat4&0&0\\0&0&2\\0&2&3\epmat$ for any $k,l\in \mathds{Z}_{\geq0}$ such that $2k-l\geq0$.
    \end{enumerate} 
    Observe that any matrix of the forms (7) to (12) are nothing but transposes of some matrix of the forms (1) to (6) and hence their corresponding entries are pre-determined by the choices above.
    \\Let us name the matrices obtained above as 
    $$A_4=\begin{pmatrix}
            0&0&2\\0&4&0\\2&0&3
        \end{pmatrix}~ A_5=\begin{pmatrix}
            0&0&2\\4&0&0\\0&2&3
        \end{pmatrix}~ A_6=\bpmat4&0&0\\0&0&2\\0&2&3\epmat$$
    \\ This allows us to infer that for any $X$ in $Im(\Psi\circ\tilde{\Phi}),$ such that it does not dominate $A_1,A_2\text{ and }A_3$, $2X\in \mathds{Z}_{\geq0}\text{linear span}\{A_4,A_5,A_6,A_5'\}$.
    \\We now record the fusion rules of these matrices which will play a pivotal role in the final computation of the quantum outer automorphism group.
    \\The fusion rules for:
    \begin{enumerate}
        \item $A_1$ : $A_1A_i=A_i$ for all i$\in \{1,...,6\}$;
    \item $A_2$ : 
 $A_2^2=A_1~~A_2A_3=A_3~~A_2A_4=A_5'~~A_2A_5=A_6~~A_2A_5'=A_4~~A_2A_6=A_5$;
    \item $A_3$ : $A_3A_2=A_3~~A_3^2=A_1+A_2+A_3~~2A_3A_4=4A_3+A_5+A_6~2A_3A_5=4A_3+A_4+A_5'~\\2A_3A_6=4A_3+A_5'+A_4~~2A_3A_5'=4A_3+A_5+A_6$;
    \item $A_4$: $A_4A_2=A_5~~2A_4A_3=4A_3+A_5+A_6~~A_4^2=4A_1+3A_4~~A_4A_5=4A_2+3A_5~\\A_4A_6=6A_3+A_5'~A_4A_5'=6A_3+A_6$;
    \item $A_5$ : $A_5A_2=A_4~~A_5^2=6A_3+A_5'~~2A_5A_3=4A_3+A_5'+A_6~~A_5A_4=6A_3+A_6~~\\A_5A_6=4A_2+3A_5~~A_5A_5'=4A_1+3A_4$;
    \item $A_6$ : $A_6A_2=A_5'~~A_6^2=4A_1+3A_6~~2A_6A_3=4A_3+A_4+A_5~~A_6A_4=6A_3+A_5~~\\A_6A_5=6A_3+A_5~A_6A_5'=4A_2+3A_5'$;
    \item $A_5'$ : $A_5'A_2=A_4~~2A_5'A_3=4A_3+A_5+A_6~~A_5'A_4=6A_3+A_6~~\\A_5'A_5=4A_1+3A6~~A
    _5'A_6=6A_3+A_4~~A_5'^2=6A_3+A_5$;
    \end{enumerate}
    Before we proceed further we take a moment to notice that if any of the matrices $A_i$ for any $i\in \{4,5,6\}$ belongs to $Im(\Psi\circ\tilde{\Phi})$ then, the fusion rules above yield that the other matrices have to be in the image as well.\\
    This can be proven very easily due to Lemma \ref{general_lemma} as follows.\\
    Without loss of generality let us assume that $A_4\in Im(\Psi\circ\tilde{\Phi})$. Since $A_1,A_2\text{ and }A_3\in Im(\Psi\circ\tilde{\Phi})$, hence $A_5=A_4A_2$ is in the image which ensures that $A_5'$ is in the image as well. Finally $A_6=A_2A_5$ which means $A_6$ is also in the image.
    \\But we can infer from the following lemmas that the fusion rules imply none of the matrices $A_i,i=4,\cdots,6$ belong to ${\rm Im}(\Psi\circ\tilde{\Phi})$.     \blmma\label{A_4,s irreducible subobjects}
     $A_4$ does not belong to ${\rm Im}(\Psi\circ\tilde{\Phi})$.
     \elmma
     \begin{proof}
     Let there exist an object $(\al,\clh)\in {\rm Rep}(\clq)$ such that $\Psi\circ\tilde{\Phi}(\al_4,\clh_{\al_4})=A_4.$
         We know from the fusion rules above that 
         \begin{align*}
             &A_4^2=4A_1+3A_4\\
             \text{or, }&\Psi\circ\tilde{\Phi}(\al_4\otimes \al_4)=\Psi\circ\tilde{\Phi}(4\alpha_1\oplus 3\alpha_4)\\
             \text{or, }&\tilde{\Phi}(\al_4\otimes \al_4)\cong\tilde{\Phi}(4\al_1\oplus 3\al_4)\\
             \text{or, }&\al_4\otimes\al_4\cong 4\al_1\oplus3\al_4.
         \end{align*}
         Thus, $\al_4$ cannot be an irreducible object in Rep($\clq)$. So, $\al_4\cong \oplus_{i=1}^rx_i$ for $x_i\in {\rm Rep}(\clq)$. Further we can infer that $\mathscr{C}_4:=\mathds{Z}_{\geq 0}\{\al_1,x_i:1\leq i\leq r\}$ forms a full $C^*$-tensor subcategory of ${\rm Rep}(\clq)$. The self-adjointness of $A_4$ implies the self-conjugacy of $\al_4$ and hence the rigidity of the sub-category $\mathscr{C}_4$ follows trivially.
         Hence using Tannaka-Krein Duality theorem on $\mathscr{C}_4$, we obtain a discrete quantum subgroup $\mathcal{S}_{\alpha_4}$ of $\clq$.
         Therefore, looking through a purely algebraic sense, we can say that $\cls_{\al_4}$ is a Hopf algebra which is a Hopf subalgebra of $\clq.$\\
         Now, let us try to look into the dimensions of the irreducible sub-objects of $\al_4$ in Rep($\clq$).
         Since the rank of $A_4=4$, hence the rank of $\tilde{\Phi}(\al_4,\clh_{\al_4})=4$, which implies the dimension of the representation $\al_4$ is also 4. Hence we have the following possibilities:
         \begin{enumerate}
             \item Let $\al_4=x_1\oplus x_2$ with dim($x_1)=1$ and dim$(x_2)=3$. This implies dim($\cls_{\al_4})=1+1+9=11$. Therefore, Theorem 3.5 of \cite{hopf_prime} implies that $\cls_{\al_4}\cong \mathds{C}[\mathds{Z}_{11}]$ thereby making $\cls_{\al_4}$ commutative. This contradicts the existence of the 3 dimensional irreducible sub-object of $\al_4$ of $\cls_{\al_4}$ which was our starting assumption.
             \item Let $\al_4=x_1\oplus x_2 $ each irreducible of dimension 2. Then we have th following sub-cases:
                \begin{enumerate}
                    \item $x_1$ and $x_2$ are dual to each other. Since by Frobenius reciprocity theorem $\al_1$ can appear in $x_1\otimes \overline{x_1}$ only once, $x_1\otimes \overline{x_1}=\al_1\oplus mx_1\oplus n x_2$ for some non negative integers m and n. Thus, calculating the dimensions of each side gives $4=1+m2+n2$ which can never happen.
                    \item $x_1$ and $x_2$ are both self conjugate. For each $1\leq i\leq 2$, $x_i\otimes x_i=\al_1\oplus mx_1\oplus nx_2$ which again yields a contradiction.                    
                \end{enumerate}
                \item Let $\al_4=x_1\oplus x_2\oplus x_3$ irreducible objects such that dim($x_3)=2 $ and the others have dimension 1 each. Then dim ($\cls_{\al_4})=1+1+1+4=7$ which again will lead to a contradiction as $\cls_{\al_4}$ will be commutative with a 2-dimensional irreducible representation.
            \end{enumerate}
            Hence, the possible dimension of any irreducible sub-object of $\al_4$ in Rep$(\clq)$ is 1. But this also means that $\Psi\circ\tilde{\Phi}(x_i)$ has rank 1 and hence is either $A_1$ or $A_2$. But $A_4$ does not dominate $A_1$ or $A_2$, which contradicts the fact that $x_i$ is a sub-object of $\al_4$. 
            Therefore, $\al_4$ cannot be an object of ${\rm Rep}(\clq)$ and hence $A_4\notin {\rm Im}(\Psi\circ\tilde{\Phi}).$
        \end{proof}
        We also have the following 
        \blmma\label{A_6's irreducible subobjects}
        $A_5,A_5'\text{ and }A_6$ are not in ${\rm Im}(\Psi\circ\tilde{\Phi})$. 
     \elmma
        \begin{proof}
            $A_6$ is also a self adjoint matrix and from the fusion rules we have $A_6^2=4A_1+3A_6.$ Therefore, we can go along the lines of the proof of the previous lemma to show that $A_6\notin{\rm Im}(\Psi\circ\tilde{\Phi})$. Again if $A_5$ is in the image then, $A_4=A_5A_2$ is also in the image which is a clear contradiction to Lemma \ref{A_4,s irreducible subobjects}. The observation that $A_5'\notin {\rm Im}(\Psi\circ\tilde{\Phi})$ is now a trivial implication.
        \end{proof}

        Our main goal now is to find the irreducible representations of $\clq$ upto isomorphism as that will enable us to infer about ${\rm Rep}(\clq)$.
       \\Before we proceed further let us fix the following notations for ease of computation:\\
       Let $G \text{ be used to denote any arbitrary object } \in \mathcal{G}:=\mathds{Z}_{\geq0}\{A_1,A_2,A_3\}$. Recall from the fusion rules that $\mathds{Z}_{\geq 0}\{A_1,A_2,A_3\}$ is closed under finite sums, finite products and the adjoint or transpose operation. Hence for any $G_1,G_2\in \clg$, we have $G_1+G_2,G_1G_2,G_2G_1,G_1'\text{ and }G_2'\in \clg.$

       \noindent Before we state the following lemma let us try to motivate it through the following observations:
       \begin{enumerate}
           \item One can easily check that $A_4+A_6=4A_1+2A_3$ and $A_5+A_5'=4A_2+2A_3.$
           \item Let $X=G_1+kA_4$ for some $G_1\in \clg$. Then we have
           \begin{align*}
               &X+A_2\{(A_2X)'\}\\
               =&G_1+kA_4+A_2(A_2G_1)'+kA_2(A_2A_4)'\\
               =&G_1+kA_4+G_2+kA_6 ~ \text{ where }G_2:=A_2(A_2G_1)'\\
               =&G_1+G_2+4kA_1+2kA_3
           \end{align*}
           Since $G_1 \text{ and }A_2\in \clg,$ hence $G_2\in \clg$, which further implies $X+A_2\{(A_2X)'\} \in \clg$.
           \item Let $X=G+kA_5$ for some $G_1\in \clg$.Then $$
                          X+X'
               =G+kA_5+G'+kA_5'
               =G+G'+k4A_2+k2A_3\in \clg. $$ Similarly if $X=G+kA_5'$, then also $X+X'\in \clg.$
           \item Let $X=G_1+kA_6$ for some $G_1\in \clg$. Then we have 
           \begin{align*}
               &X+A_2\{(A_2X)'\}\\
               =&G_1+kA_6+A_2(A_2G_1)'+kA_2(A_2A_6)'\\
               =&G_1+G_2+k4A_1+k2A_3~~\text{ where }G_2:=A_2(A_2G_1)'
           \end{align*} which implies $X+A_2\{(A_2X)'\}\in \clg.$
           \item Let $X=G_1+(2k-l)A_4+lA_5$ for some $G_1\in \clg$. Let us define $Y:=X+X'$. From (3) 
           we have $Y=G_1+G_1'+2(2k-l)A_4+l(4A_2+2A_3)=G_2+2(2k-l)A_4$ for $G_2:=G_1+G_1'+l(4A_2+2A_3)\in \clg$.
           Using (2) we have $Y+A_2\{(A_2Y)'\}\in \clg.$
           \item Let $X=G_1+(2k-l)A_5'+lA_6$ for some $G_1\in \clg$. Let us define $Y:=X+X'.$ Then we have $Y=G_1+G_1'+(2k-l)\{A_5'+A_5\}+2lA_6=G_1+G_1'+(2k-l)\{4A_2+2A_3\}+2lA_6=G_2+2lA_6$ for $G_2:=G_1+G_1'+(2k-l)\{A_5'+A_5\}+2lA_6=G_1+G_1'+(2k-l)\{4A_2+2A_3\}\in \clg$. Using (3) we have $Y+A_2\{(A_2Y)'\}\in \clg.$
           \item Let $X=G_1+(2k-l)A_4+lA_5$ for some $G_1\in \clg$.Then $X'$ is of the form seen in (5) and hence for $Y=X'+X$, we have $Y+A_2\{(A_2Y)'\}\in \clg.$
           \item Let $X=G_1+(2k-l)A_5+lA_6$ for some $G_1\in \clg$. Then $X'$ is of the form seen in (6) and hence for $Y=X'+X$, we have $Y+A_2\{(A_2Y)'\}\in \clg.$
        \end{enumerate}

         The above observations play a crucial role in the 
         \blmma \label{linear_comb_A_i_in_image}
        Let \begin{equation}\label{expr}
        \begin{split}X&=G+\sum_{i=4}^6m_iA_i+m_7A_5'+m_7\{(2k_7-l_7)A_4+l_7A_5\}+m_8\{(2k_8-l_8)A_5'+l_8A_6\}\\&+m_9\{(2k_9-l_9)A_5+l_9A_6\}+m_{10}\{(2k_{10}-l_{10})A_4+l_{10}A_5'\}\end{split}\end{equation}
        be in the image of $\Psi\circ\tilde{\Phi}$ corresponding to some object $(\al,\clh_{\al})\in {\rm Rep}(\clq)$ where $G\in \clg$. Then any irreducible sub-object of $\al\in \mathds{Z}_{\geq0}\{\al_1,\al_2,\al_3\}.$
        \elmma
        \begin{proof}
             $X=G+m_4A_4+m_5A_5+m_6A_6+m_7A_5'+m_7\{(2k_7-l_7)A_4+l_7A_5\}+m_8\{(2K_8-l_8)A_5'+l_8A_6\}+m_9\{(2k_9-l_9)A_5+l_9A_6\}+m_{10}\{(2K_{10}-l_{10})A_4+l_{10}A_5'\}$.
             \\This implies 
             \begin{align*}
             &A_2\{(A_2X)'\}=G_1+m_4A_6+m_5A_5+m_6A_4+m_7A_5'+m_7\{(2k_7-l_7)A_6+l_7A_5\}\\&~~~~~~~~~~~~~~~~~~~~~~~~~~~+m_8\{(2K_8-l_8)A_5'+l_8A_4\}+m_9\{(2k_9-l_9)A_5+l_9A_4\}\\&~~~~~~~~~~~~~~~~~~~~~~~~~~~+m_{10}\{(2K_{10}-l_{10})A_6+l_{10}A_5'\}~~~~~\text{ where }G_2\in \clg\\
             \text{or, }&X+A_2\{(A_2X)'\}=G+G_1+m_4\{A_4+A_6\}+2m_5A_5+m_6\{A_6+A_4\}+2m_7A_5'\\
             &+m_7(2k_7-l_7)\{A_4+A_6\}+2m_7l_7A_5+2m_8(2k_8-l_8)A_5'+m_8l_8\{A_6+A_4\}+2m_9(2k_9-l_9)A_5\\&+m_9l_9\{A_6+A_4\}+m_{10}(2k_{10}-l_{10})\{A_4+A_6\}+2l_{10}A_5'\\
             \text{or, }&X+A_2\{(A_2X)'\}=G_3+\{2m_5+2m_7l_7+2m_9(2k_9-l_9)\}A_5\\&~~~~~+\{2m_7+2m_8(2k_8-l_8)+2m_{10}l_{10}\}A_5'~~~~~\text{ where }G_3\in \clg\\
             \text{or, }&X+A_2\{(A_2X)'\}+X'+\{A_2X\}A_2=G_3+G_3'+\{2m_5+2m_7l_7+2m_9(2k_9-l_9)\}\{4A_1+2A_3\}\\&~~~~~~+\{2m_7+2m_8(2k_8-l_8)+2m_{10}l_{10}\}\{4A_2+2A_3\}\\
             \text{or, }&X+A_2\{(A_2X)'\}+X'+\{A_2X\}A_2\in \clg.
             \end{align*}
             So, \begin{equation}\label{eqn_1}X+A_2\{(A_2X)'\}+X'+\{A_2X\}A_2=m_1A_1+m_2A_2+m_3A_3 \text { for }m_1, m_2, m_3\in \mathds{Z}_{\geq0} \end{equation}
             We have $\al\in {\rm Rep}(\clq)$ such that $X=\Psi\circ\tilde{\Phi}(\al)$. Then  $A_2\{(A_2X)'\}\in Im(\Psi\circ\tilde{\Phi})$. Let $\bta\in {\rm Rep}(\clq)$ such that $\Psi\circ\tilde{\Phi}(\beta)=A_2\{(A_2X)'\}$.\\
             Then equation \eqref{eqn_1}  gives \begin{align*}
                 &\Psi\circ\tilde{\Phi}(\al\oplus\beta\oplus \hat{\al}\oplus\hat{\beta})=\Psi\circ\tilde{\Phi}(m_1\al_1\oplus m_2\al_2\oplus m_2\al_3)\\
                 \text{or, }&\al\oplus\beta\oplus \hat{\al}\oplus\hat{\beta}\cong m_1\tilde{\pi}(\al_1)\oplus m_2\tilde{\pi}(\al_2)\oplus m_3\tilde{\pi}(\al_3)\\
                 \text{or, }& \text{ the set of in-equivalent irreducible sub-objects of }\al\oplus\beta\oplus \hat{\al}\oplus\hat{\beta} \text{ in Rep}(\clq)\\
                 &= \{\tilde{\pi}(\al_1),\tilde{\pi}(\al_2),\tilde{\pi}(\al_3)\}\end{align*}
                Hence, any irreducible sub-object of $\al\in \mathds{Z}_{\geq 0}\{\tilde{\pi}(\al_1), \tilde{\pi}(\al_2),\tilde{\pi}(\al_3)\}.$
            Hence Proved.
    \end{proof}
    \noindent \textbf{ Proof of Theorem \ref{maximal_S3}}:
    \\Let us consider any irreducible object $(\al,\clh_{\al})\in {\rm Rep}(\clq)$. Then, $X=\Psi\circ\tilde{\Phi}(2\al)$ will be of the form \eqref{expr} given in Lemma \ref{linear_comb_A_i_in_image}. Hence, we have $\al\in \mathds{Z}_{\geq 0}\{\tilde{\pi}(\al_1), \tilde{\pi}(\al_2),\tilde{\pi}(\al_3)\}.$ We have obtained that ${\rm Rep}(\clq)\cong \mathds{Z}_{\geq0}\{\tilde{\pi}(\al_1),\tilde{\pi}(\al_2),\tilde{\pi}(\al_3)\}\cong\mathds{Z}_{\geq0}\{(\al_1),(\al_2),(\al_3)\} \cong {\rm Rep}(C^*(S_3)).$ Therefore, $\clq$ and $C^*(S_3)$ are monoidally equivalent discrete quantum groups. \qed
    \subsection{Outer action of ${\rm Rep}(\clq)$ on ${\rm Bimod}(\cla)): {\rm Out}(\cla)\subseteq {\rm Im}({\rm Rep}(\clq))$}
   Now let us consider a discrete quantum group $\clq$ with an outer action (in the sense of Definition \ref{our_action}) $\Gamma:{\rm Rep}(\clq)\to {\rm Bimod}(\cla)$ such that the ${\rm Im}(\Gamma)$ contains the image of ${\rm Out}(\cla)$. Hence, ${}_{\cla,\al_2}L^2(\cla,\tau)\otimes {\clh_{\al_2}}\in {\rm Im}(\Gamma)$, which implies the existence of an object $\beta\in {\rm Rep}(\clq)$ which serves as the pre-image of ${}_{\cla,\al_2}L^2(\cla,\tau)\otimes {\clh_{\al_2}}$ such that $A_2=\Psi\circ\Gamma(\beta)\in {\rm Im}(\Psi\circ\Gamma)$. Since the unit object $\mathds{1}\in {\rm Rep}(\clq),$ hence we already have $A_1\in {Im}(\Psi\circ\Gamma).$\\ 
   Following the same algorithm as employed previously, we can infer that for any object $X\in {\rm Im}(\Psi\circ\Gamma)$ we have $2X=m_1A_1+m_2A_2+m_3A_3+m_4A_4+m_5A_5+m_6A_6+m_7A_5'+m_7\{(2k_7-l_7)A_4+l_7A_5\}\\+m_8\{(2K_8-l_8)A_5'+l_8A_6\}+m_9\{(2k_9-l_9)A_5+l_9A_6\}+m_{10}\{(2K_{10}-l_{10})A_4+l_{10}A_5'\}.$
   \bthm\label{A3_image} 
   Let there exists an element $Z$ in $\Psi\circ\Gamma$ such that $Z\notin\mathds{Z}_{\geq0}\{A_1,A_2\}.$ Then we have ${\rm Im}(\Psi\circ\Gamma)=\mathds{Z}_{\geq0}\{A_1,A_2,A_3\}.$ 
   \ethm
   \begin{proof}
    Let $Z$ be any arbitrary object in ${\rm Im}(\Psi\circ\Gamma)$. Then, $2Z$ will be of the form \eqref{expr} as in Lemma \ref{linear_comb_A_i_in_image}. Therefore, 
    \begin{align*}
        2Z+A_2\{(A_22Z)'\}+2Z'+\{A_22Z\}A_2=m_1A_1+m_2A_2+m_3A_3\text{ where}m_1,m_2\geq0\text{ and }m_3>0.
    \end{align*}
    Therefore, there exists an irreducible sub-object of the pre-image of $Z$, namely, $\al$ such that $\tilde{Z}=\Psi\circ\Gamma(\al)=m'_1A_1+m'_2A_2+m'_3A_3$ with $m'_3>0.$ Then, because of the self-adjointness of $\Psi\circ \Gamma(\al)$, we have \begin{align*}
        &\Psi\circ\Gamma(\al\otimes \overline{\al})={m'_1}^2A_1+{m'_2}^2A_1+{m'_3}^2\{A_1+A_2+A_3\}+2m'_1m'_2A_2+2m'_2m'_3A_3+2m'_1m'_3A_3\\
        \text{or, }&\Psi\circ\Gamma(\al\ot\overline{\al})+m'_1m'_3A_1+m'_2m'_3A_2+2{m'_1}^2A_1+2m'_1m'_2A_1+2{m'_2}^2A_2\\
        &={m'_1}^2A_1+{m'_2}^2A_1+{m'_3}^2A_1+{m'_3}^2A_2+\{m'_1m'_3A_1+m'_2m'_3A_2+{m'_3}^2A_3\}+\\&\{2{m'_1}^2A_1+2m'_1m'_2A_2+2m'_1m'_3A_3\}+\{2m'_1m'_2A_1+2{m'_2}^2A_2+2m'_2m'_3A_3\}\\
\text{or, }&\Psi\circ\Gamma(\al\ot\overline{\al}\oplus \{m'_1m'_3+2{m'_1}^2+2m'_1m'_2\}{\mathds{1}}\oplus\{m'_2m'_3+2m'_1m'_2+2{m'_2}^2\}\bta)=\\
        &\Psi\circ\Gamma(\{{m'_1}^2+{m'_2}^2+{m'_3}^2\}{\mathds{1}}\oplus \{+{m'_3}^2+2m'_1m'_2\}\bta\oplus m'_3\al\oplus 2m'_1\al\oplus 2m'_2\al)\\
        \text{or, }&\al\ot\overline{\al}\oplus \{m'_1m'_3+2{m'_1}^2+2m'_1m'_2\}
        {\rm 
        Id}\oplus\{m'_2m'_3+2m'_1m'_2+2{m'_2}^2\}\bta\cong\\&\{{m'_1}^2+{m'_3}^2\}{\mathds{1}}\oplus \{{m'_2}^2+{m'_3}^2+2m'_1m'_2\}\bta\oplus m'_3\al\oplus 2m'_1\al\oplus 2m'_2\al\\
        \text{or, }&{\rm dim(Mor)}({\mathds{1}},\al\ot\overline{\al}\oplus \{m'_1m'_3+2{m'_1}^2+2m'_1m'_2\}
        {\rm 
        Id}\oplus\{m'_2m'_3+2m'_1m'_2+2{m'_2}^2\}\bta)={\rm dim(Mor)}({\mathds{1}},\\&\{{m'_1}^2+{m'_2}^2+{m'_3}^2\}{\mathds{1}}\oplus \{{m'_3}^2+2m'_1m'_2\}\bta\oplus m'_3\al\oplus 2m'_1\al\oplus 2m'_2\al)\\
        \text{or, }&1+m'_1m'_3+2{m'_1}^2+2m'_1m'_2={m'_1}^2+{m'_2}^2+{m'_3}^2.\\
    \end{align*}
    Therefore, we have \begin{equation}\label{star}
        1={m'_3}^2-{m'_1}^2+{m'_2}^2-2m'_1m'_2-m'_1m'_3
    \end{equation}
    Now, observe that
    $$
        \Psi\circ\Gamma(\bta\ot\bta)=A_2^2=A_1=\Psi\circ\Gamma({\mathds{1}})
        \text{or, }\bta\ot\bta\cong{\mathds{1}}
        \text{or, }\bta\otimes \bar{\bta}\cong{\mathds{1}}.
    $$ This implies, if $\al$ is irreducible, then $\al\ot\bta$ is also irreducible. Now, $\Psi\circ\Gamma(\al\ot\bta)=A_2\tilde{Z}=m'_2A_1+m'_1A_2+m'_3A_3$. Therefore, equation \eqref{star} gives us \begin{equation}\label{star2}
        1={m'_3}^2-{m'_2}^2+{m'_1}^2-2m'_1m'_2-m'_2m'_3.
    \end{equation} Putting, both equations together we have 
    \begin{align*}
        &{m'_1}^2 -{m'_2}^2+m'_1m'_3={m'_2}^2-{m'_1}^2+m'_2m'_3\\
        \text{or, }&2(m'_1+m'_2)(m'_1-m'_2)+m'_3(m'_1-m'_2)=0\\
        \text{or, }&(m'_1-m'_2)(2m'_1+2m'_2+m'_3)=0\text{ which implies }
        m'_1=m'_2
    \end{align*}
    Therefore, if $\al$ is irreducible then $\tilde{Z}=m'_1A_1+m'_2A_2+m'_3A_3$ with $m'_1=m'_2$ such that equation \eqref{star} holds. Thus, \begin{align*}
        &1={m'_3}^2-{m'_2}^2+{m'_2}^2-2m'_2m'_2-m'_2m'_3\\
        \text{or, }&1={m'_3}^2-2{m'_2}^2-m'_2m'_3 \end{align*}
    One can check easily that the only non negative integral values of $(m'_2,m'_3)$ which satisfies this equation is $(0,1)$. This gives us that ${\rm Im}(\Psi\circ\Gamma)(\al)=\tilde{Z}=A_3.$
    \end{proof}
    Thus, we have the following 
    \bcrlre\label{contains_Out}
    Let $\clq$ be any DQG with an action $\Gamma:{\rm Rep}(\clq)\to {\rm Bimod}(\cla)$ such that ${\rm Im}(\Gamma)$ contains ${\rm Out}(\cla)$. Then $\clq$ is monoidally equivalent to ${\rm Out}(\cla)$ or $C^*(S_3)$. \ecrlre
    \begin{proof}
        Let us first assume that ${\rm Im}(\Gamma)=\mathds{Z}_{\geq0}\{A_1,A_2\}$. Then the injectivity of the functor $\Gamma$ implies that ${\rm Rep}(\clq)$ is monoidally equivalent to ${\rm Rep(Out}(\cla))$. On other other hand if the inclusion is strict, then Theorem \ref{A3_image} implies that ${\rm Rep}(\clq)$ is monoidally equivalent to ${\rm Rep}(C^*(S_3))$. Thus, $\clq$ is monoidally equivalent to ${\rm Out}(\cla)$ or $C^*(S_3)$.  
    \end{proof} We end with a brief discussion on a particular choice of the DQG $\clq$ considered in the results obtained so far. 
    \bdfn[Definition 1.1, \cite{etingof_gelaki}]
    Two groups $G_1$ and $G_2$ are called isocategorical if ${\rm Rep}(G_1)$ is equivalent to ${\rm Rep}(G_2)$ as a tensor category.
    \edfn Moreover, we call a finite group $G$ to be categorically rigid if any group isocategorical to $G$ is actually isomorphic to $G$. Corollary 1.4 of \cite{etingof_gelaki} implies that $S_3$ is a categorically rigid group. 
   Therefore, we have the following 
    \bthm\label{grp_algebra}
    Let $\clq=C^*(H)$ for a finite group $H$, such that ${\rm Rep}(\clq)$ satisfies Theorem \ref{maximal_S3} and Corollary \ref{contains_Out}. Then, $H$ is isomorphic to $S_3$ and hence $\clq\cong C^*(S_3)$.
    \ethm       \bibliographystyle{amsalpha} 
\bibliography{reference} 

\end{document}